# A Monadic, Functional Implementation of Real Numbers[*]

Russell O'Connor

Institute for Computing and Information Science
Faculty of Science
Radboud University Nijmegen

*Email:* `r.oconnor@cs.ru.nl`

Abstract. Large scale real number computation is an essential ingredient in several modern mathematical proofs. Because such lengthy computations cannot be verified by hand, some mathematicians want to use software proof assistants to verify the correctness of these proofs. This paper develops a new implementation of the constructive real numbers and elementary functions for such proofs by using the monad properties of the completion operation on metric spaces. Bishop and Bridges's [1] notion of regular sequences is generalized to, what I call, regular functions which form the completion of any metric space. Using the monad operations, continuous functions on length spaces (a common subclass of metric spaces) are created by lifting continuous functions on the original space. A prototype Haskell implementation has been created. I believe that this approach yields a real number library that is reasonably efficient for computation, and still simple enough to easily verify its correctness.

## 0. Licence

## 1. Introduction

Several mathematical theorems rely on numerical computation of real number values in their proofs. One example is the disproof of the Mertens conjecture [15]. The Mertens conjecture claims that the absolute values of the partial sums of the Möbius function are bounded by the square root function. The conjecture implies the Riemann hypothesis. Odlyzko and te Riele's disproof involves computing the first hundred decimal digits of the first two thousand zeros of the Riemann zeta function.

Another example is Hales's famous proof [6] of the Kepler conjecture. This proof involves verifying thousands of nonlinear inequalities over the real numbers. These inequalities were verified by special purpose software. Unfortunately, the peer review process has failed to fully accept the proof as correct. After four years, the twelve referees could only claim to be 99% certain of the correctness of the proof. In response, Hales has created the *Flyspeck* project, a project to create a software verified proof of the Kepler conjecture. A full software verified proof will verify both the mathematical and computational parts of the proof together.

Many other proofs rely on real number computation. All of these proofs require a software verified implementation of real number arithmetic before these proofs can be verified by computer. By implementation, I mean a software verified library of continuous functions on $\mathbb{R}$ with which computations can be done to arbitrary precision. One such library has been created by Cruz-Filipe as a product of his software verified constructive proof of the fundamental theorem of calculus [4].

When constructive mathematics is used for reasoning, algorithms are contained inside the proofs. These algorithms can be executed in a functional language. One can see this as part of the Bishop's program to see constructive mathematics as a programming language. Because Cruz-Filipe's proof is constructive, the functions in his library can actually be evaluated. However, Cruz-Filipe's construction was not designed with efficiency in mind, so, unfortunately, evaluation is not practical [5]. Therefore, if one wants to be able to handle the problems encountered in the aforementioned mathematical proofs, a new library is needed.

This paper proposes a new constructive implementation of real numbers. An implementation that attempts to be simple and elegant, and at the same time practical enough to run the computations needed by the proof of the Kepler conjecture and other mathematical theorems. This work contributes not only to the *Flyspeck* project, but to any project requiring a verified library of functions for exact real arithmetic.

Section 2 introduces the notion of a regular function. Regular functions are a generalization of Bishop and Bridges's notion of a regular sequence [1]. Section 3 discusses how regular functions, under a suitable equivalence relation, form a completion of any metric space.

Section 4 shows that the completion operation of a metric space forms a monad. A monad is an abstract structure consisting of a type constructor and three functions: unit, map, and join. These three functions must satisfy seven monad laws, which are given in section 4.1. Monadic programming techniques have been shown to be useful for programming tasks such as I/O, exceptions, non-determinism, and many more [11][17]. The fact that the completion operation forms a monad is not so exciting, but what is exciting is that these monad functions are very useful for defining functions over completed metric spaces by lifting uniformly continuous functions on the original metric space. Section 5 shows how this process can be used to define the real numbers and define the elementary functions over $\mathbb{R}$.

Because one goal of this work is to be able to evaluate expressions to a given precision in a reasonable amount of time, section 6 gives some minor enhancements to the implementation to improve performance. Section 7 compares the run-time of a Haskell prototype implementation to other systems.

A second goal of this work is to eventually verify the proofs presented in this paper using a software proof assistant. To this end, I have included many theorems with detailed proofs in order to illustrate that verifying these proofs by software will be relatively easy. Those readers who are only interested in the data structures and operations may wish to skip over the details of the proofs.

I use constructive mathematics for all the proofs in this paper; therefore, I remain agnostic on whether or not functions are assumed to be computable. Of course, the goal is to make an implementation of the functions defined in this paper so the computational interpretation is worth keeping in mind; however both the constructive and classical interpretations are sound.

## 1.1. Notation.

The strictly positive rationals are denoted by $\mathbb{Q}^+$, while the non-negative rationals are denoted by $\mathbb{Q}^{0+}$; similar notation is used for reals and natural numbers. Equivalence relations are denoted by $\asymp$, while $=$ is reserved for intensional equality. I will be using sequences and functions to represent real numbers. It is possible for two different sequences or functions to represent the same real number, so it is important to distinguish between equality as functions and equivalence as real numbers.

The type of propositions is denoted by $\star$, and relations by functions to $\star$. If reasoning classically, $\star$ may be considered as a Boolean type.

In functional programming, functions are typically curried. Therefore, for a binary operator that would normally have the type $X \times X \Rightarrow X$, I will instead use the type $X \Rightarrow X \Rightarrow X$; however, I will sometimes write function application as $f(x, y)$ for clarity.

As is also typical in functional programming, I will write anonymous functions using lambda expressions. For example, I will write the reciprocal function as $\lambda x . x^{-1}$.

## 2. Real Numbers as Regular Functions of Rationals

Abstractly, real numbers can be defined as any complete, Archimedean, ordered field; however, in order to show that such a structure exists, it is necessary to produce a model of the real numbers. One common model of the real numbers consists of equivalence classes of Cauchy sequences of rational numbers.

Bishop and Bridges [1] define the real numbers to be regular sequences of rational numbers with an equivalence relation. A sequence $x_n$ is *regular* if

$$\forall n\, m : \mathbb{N}^+, |x_m - x_n| \leq m^{-1} + n^{-1},$$

and two regular sequences $x$ and $y$ are equivalent $(x \asymp y)$ if

$$\forall n \colon \mathrm{N}^+, |x_n - y_n| \le 2\, n^{-1}.$$

Regular sequences, under this equivalence relation, are isomorphic to Cauchy sequences. Every regular sequence is a Cauchy sequence, and every Cauchy sequence can be transformed into a regular sequence by using its modulus of convergence to select a suitable subsequence.

One can think of a regular sequence as a function that approximates a real number by taking a positive number $n$ to a rational number $x_n$ that is within $n^{-1}$ of the real number the sequence represents. Instead of using such coarse grained approximations, one can generalize the concept of regular sequences to regular functions. A function $x \colon \mathbb{Q}^+ \Rightarrow \mathbb{Q}$ is *regular* if

$$\forall \varepsilon_1\, \varepsilon_2, |x(\varepsilon_1) - x(\varepsilon_2)| \le \varepsilon_1 + \varepsilon_2,$$

and two regular functions $x$ and $y$ are equivalent if

$$\forall \varepsilon, |x(\varepsilon) - y(\varepsilon)| \le 2\, \varepsilon.$$

Regular functions are isomorphic to regular sequences. Any regular function can be transformed into an regular sequence by considering $\lambda n \colon \mathrm{N}^+.x(n^{-1})$, and any regular sequence can be transformed into a regular function, that maps $\varepsilon$ to some $x_n$ such that $n^{-1} \le \varepsilon$.

One can think of a regular function as a function that approximates a real number by taking a positive rational number $\varepsilon$ to a rational number that is within $\varepsilon$ of the real number the function represents. Regular functions allow more fine grained approximations than regular sequences allow. This fine granularity can prevent unnecessary over approximation when doing calculations. An example of this occurs in the implementation of multiplication in section 5.

Regular functions can be used not only to construct the reals, but they can be used to complete any metric space.

## 3. Metric Spaces

A typical definition of a metric space is space $X$ with a distance function $d \colon X \Rightarrow X \Rightarrow \mathbb{R}^{0+} \cup \{\infty\}$ satisfying certain properties [3]; however, the purpose of this paper is to define $\mathbb{R}$ as the completion of a metric space. This task requires a definition of metric space that does not presuppose the existence of $\mathbb{R}$. Given a traditional metric space, one can define a distance relation $B_\varepsilon(x, y) \overset{def}{=} d(x, y) \le \varepsilon$ where $\varepsilon \colon \mathbb{Q}^+$. This distance relation does not depend on $\mathbb{R}$, but it still completely characterizes the metric space.

A tuple $(X, \asymp, B)$ is a *metric space*, where $\asymp$ is an equivalence relation, $B \colon \mathbb{Q}^+ \Rightarrow X \Rightarrow X \Rightarrow \star$ is a relation respecting the equivalence relation on $X$, and the following properties hold:

1. For all $\varepsilon$, $B_\varepsilon$ is reflexive.

2. For all $\varepsilon$, $B_\varepsilon$ is symmetric.

3. For all $\varepsilon_1$, $\varepsilon_2$, $a$, $b$, and $c$, if $B_{\varepsilon_1}(a, b)$ and $B_{\varepsilon_2}(b, c)$ hold, then $B_{\varepsilon_1 + \varepsilon_2}(a, c)$ holds.

4. For all $\varepsilon$, $a$, and $b$, if $\forall \delta, B_{\varepsilon + \delta}(a, b)$ holds, then $B_\varepsilon(a, b)$ holds.

5. For all $a$ and $b$, if $\forall \varepsilon, B_\varepsilon(a, b)$ holds, then $a \asymp b$.

Given a metric space, the traditional metric can be recovered (classically) by defining $d(a, b) \overset{def}{=} \inf \{\varepsilon \mid B_\varepsilon(a, b)\}$—the infimum of the empty set is taken to be $\infty$. One can easily show that this definition of a metric space is equivalent to the traditional definition.

I use the symbol $B$ for the distance relation because if $B_\varepsilon(a)$ is considered as a predicate, then the set of points $b$ satisfying $B_\varepsilon(a)$ are exactly the points that are within $\varepsilon$ of $a$. Therefore, $B_\varepsilon(a)$ can be understood as a ball of radius $\varepsilon$ around $a$.

Property 4 requires that $B_\varepsilon(a)$ be a closed ball. While the choice between open and closed balls is not very important from a classical perspective, there are some advantages to using closed balls when reasoning constructively. Because being closed is a negative statement, it does not contain any constructive information. One can avoid unnecessary computation and allow reasoning by contradiction by using negative statements wherever possible [5].

It is useful to note that if two points are within $\varepsilon_1$ of each other, then they are also within $\varepsilon_2$ of each other whenever $\varepsilon_1 < \varepsilon_2$.

**Lemma 1.** *For any $\varepsilon$, $\delta$, $a$, and $b$, if $B_\varepsilon(a,b)$ holds, then $B_{\varepsilon+\delta}(a,b)$ holds.*

**Proof.** $B_\delta(b,b)$ holds by reflexivity. Therefore, $B_{\varepsilon+\delta}(a,b)$ holds.          □

### 3.1. Prelength Spaces.

In a metric space, just because two points $a$ and $b$ are within $l$ of each other does not mean that there exist curves of length approaching $l$ that connect these points. A metric space that has such curves is known as a length space [3]. Since the notion of curves presupposes the real numbers, I introduce the weaker notion of a prelength space. A *prelength space* is a metric space such that

$$\forall a\,b\,\varepsilon\,\delta_1\,\delta_2, \varepsilon < \delta_1 + \delta_2 \Rightarrow B_\varepsilon(a,b) \Rightarrow \exists c, B_{\delta_1}(a,c) \wedge B_{\delta_2}(c,b).$$

This property implies that if two points $a$ and $b$ are within $l$ of each other, then there is a path of points connecting $a$ to $b$ with consecutive points within $\varepsilon$ of each other whose total length is less than $l + \varepsilon$.

**Lemma 2.** *In a prelength space for any $a$, $b$, $\varepsilon$, $\delta_0$, ..., $\delta_{n-1}$ such that $B_\varepsilon(a,b)$ holds and $\varepsilon < \delta_0 + ... + \delta_{n-1}$, there exist $c_0, ..., c_n$ such that $c_0 = a$, $c_n = b$, and $\forall i < n, B_{\delta_i}(c_i, c_{i+1})$ holds.*

**Proof.** This is done by induction on $n$. The case $n = 0$ is trivial. In the other case, let $\gamma \overset{\text{def}}{=} \delta_0 + ... + \delta_{n-1} - \varepsilon$. There exists some $c_{n-1}$ such that $B_{\delta_{n-1}}(c_{n-1}, b)$ and $B_{\varepsilon - \delta_{n-1} + \frac{\gamma}{2}}(a, c_{n-1})$. By induction, there exists $c_0, ..., c_{n-2}$ such that $\forall i < n-1, B_{\delta_i}(c_i, c_{i+1})$ because $\varepsilon - \delta_{n-1} + \frac{\gamma}{2} < \delta_0 + ... + \delta_{n-2}$.          □

Prelength spaces are quite common; all the metric spaces considered in this paper are prelength spaces.

**Theorem 3.** $(\mathbb{Q}, =, B^{\mathbb{Q}})$ *is a prelength space where* $B_\varepsilon^{\mathbb{Q}}(a,b) \overset{\text{def}}{=} |a - b| \leq \varepsilon$.

**Theorem 4.** *For any closed rational interval $[a,b]$, $([a,b], =, B^{\mathbb{Q}})$ is a prelength space.*

### 3.2. Completion of a Metric Space using Regular Functions.

One can define regular functions over any metric space $X$ by using the distance relation. A function $x: \mathbb{Q}^+ \Rightarrow X$ is *regular* if

$$\forall \varepsilon_1\,\varepsilon_2, B_{\varepsilon_1 + \varepsilon_2}(x(\varepsilon_1), x(\varepsilon_2)),$$

and two regular functions $x$ and $y$ are equivalent $(x \asymp y)$ if

$$\forall \varepsilon, B_{2\varepsilon}(x(\varepsilon), y(\varepsilon)).$$

Define the distance relation on regular functions as

$$B'_\varepsilon(x,y) \overset{\text{def}}{=} \forall \delta_1\,\delta_2, B_{\varepsilon + \delta_1 + \delta_2}(x(\delta_1), y(\delta_2)).$$

Call the type of regular functions over $X$ the *completion* of $X$ or $\mathfrak{C}(X)$. The following theorems show that the completion of a metric space is a metric space.

**Theorem 5.** *If $(X, \asymp, B)$ is a metric space, then $(\mathfrak{C}(X), \asymp, B')$ is a metric space.*

**Proof.**

1. Consider arbitrary $\delta_1$ and $\delta_2$. By Lemma 1, $B_{\varepsilon + \delta_1 + \delta_2}(x(\delta_1), x(\delta_2))$ holds. Therefore, $B'_\varepsilon$ is reflexive.

2. $B'$ is symmetric because $B$ is symmetric.

3. Suppose $B'_{\varepsilon_1}(x, y)$ and $B'_{\varepsilon_2}(y, z)$ hold. Let $\gamma$ be arbitrary. Then $B_{\varepsilon_1 + \delta_1 + \frac{\gamma}{2}}(x(\delta_1), y(\frac{\gamma}{2}))$ and $B_{\varepsilon_2 + \frac{\gamma}{2} + \delta_2}(y(\frac{\gamma}{2}), z(\delta_2))$ hold. Therefore, $B_{\varepsilon_1 + \varepsilon_2 + \delta_1 + \delta_2 + \gamma}(x(\delta_1), y(\delta_2))$ holds. Since $\gamma$ was arbitrary, $B'_{\varepsilon_1 + \varepsilon_2}(x, y)$ holds.

4. Given $\varepsilon$, $x$, and $y$, suppose $\forall \delta_0, B'_{\varepsilon + \delta_0}(x, y)$ hold. Therefore, $\forall \delta_0 \, \delta_1 \, \delta_2, B_{\varepsilon + \delta_0 + \delta_1 + \delta_2}(x(\delta_1), y(\delta_2))$ holds. $\forall \delta_1 \, \delta_2, B_{\varepsilon + \delta_1 + \delta_2}(x(\delta_1), y(\delta_2))$ holds; therefore, $B'_{\varepsilon}(x, y)$ holds.

5. Given $x$ and $y$, suppose $\forall \varepsilon, B'_{\varepsilon}(x, y)$ hold. For any $\varepsilon$ and $\delta$, $B_{\varepsilon + 2\delta}(x(\delta), y(\delta))$ holds. Then $B_{2\delta}(x(\delta), y(\delta))$ holds. Therefore, $x \asymp y$.

$\square$

The following are useful lemmas about the distance function on regular functions.

**Lemma 6.** *If $B_{\varepsilon_0}(x(\varepsilon_1), y(\varepsilon_2))$ holds, then $B'_{\varepsilon_0 + \varepsilon_1 + \varepsilon_2}(x, y)$.*

**Proof.** Let $\delta_1$ and $\delta_2$ be arbitrary. $B_{\varepsilon_1 + \delta_1}(x(\delta_1), x(\varepsilon_1))$ and $B_{\varepsilon_2 + \delta_2}(y(\varepsilon_2), y(\delta_2))$ hold. Therefore, $B_{\varepsilon_0 + \varepsilon_1 + \varepsilon_2 + \delta_1 + \delta_2}(x(\delta_1), y(\delta_2))$ holds. $\square$

**Corollary 7.** *If $B_{\varepsilon_0}(x(\varepsilon_1), b)$ holds, then $B'_{\varepsilon_0 + \varepsilon_1}(x, b')$ holds where $b'$ is the regular function $b' \overset{def}{=} \lambda \varepsilon. b$.*

**Proof.** By Lemma 6 $B'_{\varepsilon_0 + \varepsilon_1 + \varepsilon_2}(x, b')$ holds for every $\varepsilon_2$. Therefore, $B'_{\varepsilon_0 + \varepsilon_1}(x, b')$ holds. $\square$

The following theorem shows that the distance relation is equivalent to $\forall \delta, B_{\varepsilon + 2\delta}(x(\delta), y(\delta))$.

**Theorem 8.** *For any $\varepsilon$, if $\forall \delta, B_{\varepsilon + 2\delta}(x(\delta), y(\delta))$ holds, then $B'_{\varepsilon}(x, y)$ holds.*

**Proof.** For any $\delta$, $B_{\varepsilon + \frac{\delta}{2}}(x(\frac{\delta}{4}), y(\frac{\delta}{4}))$ holds. By Lemma 6, $B'_{\varepsilon + \delta}(x, y)$ holds. Since $\delta$ was arbitrary, $B'_{\varepsilon}(x, y)$ holds. $\square$

The following theorem shows that the completion of a prelength space is a prelength space.

**Theorem 9.** *If $(X, \asymp, B)$ is a prelength space, then is a prelength space.*

**Proof.** By Theorem 5, all that remains is to show that the prelength space property is preserved. Suppose $\varepsilon < \delta_1 + \delta_2$, and $B'_{\varepsilon}(x, y)$ holds. Let $\gamma \overset{def}{=} \frac{\delta_1 + \delta_2 - \varepsilon}{5}$. $B_{\varepsilon + 2\gamma}(x(\gamma), y(\gamma))$ holds. There is some $c$ such that $B_{\delta_1 - \gamma}(x(\gamma), c)$ and $B_{\delta_2 - \gamma}(c, y(\gamma))$ hold. So $B'_{\delta_1}(x, c')$ and $B'_{\delta_2}(c', y)$ holds by Corollary 7 where $c' \overset{def}{=} \lambda \varepsilon. c$. $\square$

Typically the completion of a space $X$ is denoted by $\bar{X}$. Because I will wish to emphasize the monadic property of the completion, the Fraktur $\mathfrak{C}$ is used in this paper.

### 3.3. Uniform Continuity.

Given two metric spaces $(X, \asymp, B^X)$ and $(Y, \asymp, B^Y)$, a function $f \colon X \Rightarrow Y$ is *uniformly continuous* with *modulus* $\mu \colon \mathbb{Q}^+ \Rightarrow \mathbb{Q}^+$ if

$$\forall \varepsilon \, x_1 \, x_2 \colon X . B^X_{\mu(\varepsilon)}(x_1, x_2) \Rightarrow B^Y_{\varepsilon}(f(x_1), f(x_2)).$$

A function $f \colon X \Rightarrow Y$ is *uniformly continuous* if there is some $\mu$ such that $f$ is uniformly continuous with modulus $\mu$. When a function $f$ is uniformly continuous, I will write $f \colon X \to Y$ using a single-bar arrow. This indicates that $f$ is not just a function, but it is a pair $(f, \mu_f) \colon (X \Rightarrow Y) \times (\mathbb{Q}^+ \Rightarrow \mathbb{Q}^+)$ consisting of a function and its modulus of continuity. This structure is a morphism in the category of metric spaces with uniformly continuous functions between them. When $f \colon X \to Y$, I will leave the projection function implicit and write $f(x)$ instead of $\pi_1(f)(x)$.

This is the usual definition of uniform continuity where the relation between $\varepsilon$ and $\delta$ is explicitly given by $\mu$.

**Theorem 10.** *The identity function* id *is uniformly continuous with modulus* id, *and if $f \colon X \to Y$ and $g \colon Y \to Z$ are both uniformly continuous, then the composition $g \circ f \colon X \to Z$ is uniformly continuous with modulus $\mu_f \circ \mu_g$.*

**Proof.** Easy.                                                                                         □

### 3.4. Remarks.

The definition of a metric space usually requires the distance between two points always to be finite. Because I want consider the space of all uniformly continuous functions between two metric spaces as a metric space in section 4.2, I use the more liberal definition. It is worth noting that it is only the uniform space structure of a metric space that is important for the purposes of this paper. One can transform a metric space with infinite distances to one with only finite distances by defining a new metric $d'(a, b) \overset{\text{def}}{=} \min \{d(a, b), 1\}$. The uniform structure of the space does not change by using this new metric.

A prelength space is a weaker notion than a length space. For example, the rationals are not a length space because no two points can be connected by a continuous curve; however, the rationals are a prelength space. The name *prelength space* is chosen because of an analogy with precompact spaces. When a precompact space is completed it becomes a compact space; when a prelength space is completed it becomes a length space.

An alternative definition for a regular function $x \colon \mathbb{Q}^+ \Rightarrow X$ would be

$$\forall \varepsilon_1 \varepsilon_2, B_{\max(\varepsilon_1, \varepsilon_2)}(x(\varepsilon_1), x(\varepsilon_2)).$$

This alternative definition would remove the need for restricting ourselves to prelength spaces (in particular, its use in Theorem 16). However, for efficiency reasons discussed in section 6.1, I stick with the definition based on Bishop and Bridges's regular sequences.

## 4. Completion Is a Monad

Any completion constructor forms a monad in the category of prelength spaces and uniformly continuous functions between them. The following defines the monad structure for the definition of completion given in section 3.2. The unit $\colon X \to \mathfrak{C}(X)$ function is the obvious injection of $X$ into $\mathfrak{C}(X)$.

$$\text{unit}(a) \overset{\text{def}}{=} \lambda \varepsilon. a.$$

**Theorem 11.** *For any $a \colon X$,* unit$(a)$ *is a regular function, and* unit *is uniformly continuous with modulus* id.

The following two theorems show that the geometry of the space is preserved under this injection and that it is sound to think of $x(\varepsilon)$ as an approximation of $x$ within $\varepsilon$.

**Theorem 12.** *For any $a$, $b$, and $\varepsilon$, $B_\varepsilon(a, b)$ holds if and only if $B'_\varepsilon(\text{unit}(a), \text{unit}(b))$ holds.*

**Theorem 13.** *For any $\varepsilon$, $B'_\varepsilon(x, \text{unit}(x(\varepsilon)))$ holds.*

**Proof.** Let $\delta$ be arbitrary. $B_{\varepsilon + \delta}(x(\delta), x(\varepsilon))$ holds because $x$ is regular, but this is the same as $B_{\varepsilon + \delta}(x(\delta), \text{unit}(x(\varepsilon))(\delta))$. By Theorem 8, $B'_\varepsilon(x, \text{unit}(x(\varepsilon)))$ holds.                    □

The function join $\colon \mathfrak{C}(\mathfrak{C}(X)) \to \mathfrak{C}(X)$ arises from the proof that completing a space twice is the same space as completing it only once.

$$\text{join}(x) \overset{\text{def}}{=} \lambda \varepsilon. x(\tfrac{\varepsilon}{2})(\tfrac{\varepsilon}{2}).$$

**Theorem 14.** *For any $x \colon \mathfrak{C}(\mathfrak{C}(X))$,* join$(x)$ *is a regular function, and* join *is uniformly continuous with modulus* id.

**Proof.** Let $\varepsilon_1$ and $\varepsilon_2$ be arbitrary. $B_{\frac{\varepsilon_1}{2} + \frac{\varepsilon_2}{2}}(x(\frac{\varepsilon_1}{2}), x(\frac{\varepsilon_2}{2}))$ holds because $x$ is a regular function. By Theorem 13, $B_{\frac{\varepsilon_i}{2}}(x(\frac{\varepsilon_i}{2}), \text{unit}(x(\frac{\varepsilon_i}{2})(\frac{\varepsilon_i}{2})))$ hold, so all together $B_{\varepsilon_1 + \varepsilon_2}(\text{unit}(\text{join}(x)(\varepsilon_1)), \text{unit}(\text{join}(x)(\varepsilon_2)))$ holds. By Theorem 12, $B_{\varepsilon_1 + \varepsilon_2}(\text{join}(x)(\varepsilon_1), \text{join}(x)(\varepsilon_2))$ holds, so join$(x)$ is a regular function.

Now let $\varepsilon$ and $\delta$ be arbitrary, and suppose $x_1$ and $x_2$ are such that $B_\varepsilon(x_1, x_2)$ holds. $B_\delta(x_i, \mathrm{unit}(\mathrm{unit}(x_i(\frac{\delta}{2})(\frac{\delta}{2}))))$ holds by two applications of Theorem 13. Therefore, $B_{\varepsilon+2\delta}(\mathrm{unit}(\mathrm{unit}(\mathrm{join}(x_1)(\delta))), \mathrm{unit}(\mathrm{unit}(\mathrm{join}(x_2)(\delta))))$ holds. By two applications of Theorem 12, $B_{\varepsilon+2\delta}(\mathrm{join}(x_1)(\delta), \mathrm{join}(x_2)(\delta))$ holds. By Theorem 8, $B_\varepsilon(\mathrm{join}(x_1), \mathrm{join}(x_2))$ holds as required.  □

The function map : $(X \to Y) \Rightarrow (\mathfrak{C}(X) \to \mathfrak{C}(Y))$ lifts uniformly continuous functions on the base spaces to uniformly continuous on the completed spaces.

$$\mathrm{map}(f) \stackrel{def}{=} \lambda x.\, f \circ x \circ \mu_f.$$

For metric spaces, $\mathrm{map}(f)(x)$ is not always a regular function. Let $X$ be the space $\{-1\} \cup \{n^{-1} \,|\, n \colon \mathbb{N}^+\}$ and $Y$ be the space $\{-2\} \cup \{n^{-1} \,|\, n \colon \mathbb{N}^+\}$ with the induced metrics. The function $f \colon X \to Y$ sending $-1$ to $-2$ and $n^{-1}$ to $n^{-1}$ is uniformly continuous with modulus $\lambda\varepsilon.\min(\varepsilon, 1)$. Consider the regular function $x$ that sends $\varepsilon$ to $-1$ when $1 \le \varepsilon$ and sends $\varepsilon$ to some $n^{-1} \le \varepsilon$ when $\varepsilon < 1$. In this case, $\mathrm{map}(f)(x)$ is not a regular function.

Fortunately $\mathrm{map}(f)(x)$ is always regular function when $X$ is a prelength space.

**Lemma 15.** *Given $f \colon X \to Y$ where $X$ is a prelength space and $Y$ is a metric space, consider arbitrary $a$, $b$, $\varepsilon_0$, ..., $\varepsilon_{n-1}$. Let $\delta_i \stackrel{def}{=} \mu_f(\varepsilon_i)$, $\delta \stackrel{def}{=} \delta_0 + ... + \delta_{n-1}$, and $\varepsilon \stackrel{def}{=} \varepsilon_0 + ... + \varepsilon_{n-1}$. If $B_\delta(a, b)$ holds, then $B_\varepsilon(f(a), f(b))$ holds.*

**Proof.** Let $\varepsilon_n$ be arbitrary, and let $\delta_n \stackrel{def}{=} \mu_f(\varepsilon_n)$. By Lemma 2, there are $c_0$, ..., $c_{n+1}$ such that $B_{\delta_i}(c_i, c_{i+1})$ hold, $a = c_0$, and $b = c_{n+1}$. By uniform continuity $B_{\varepsilon_i}(f(c_i), f(c_{i+1}))$ hold. Thus $B_{\varepsilon+\varepsilon_n}(f(a), f(b))$ holds. Since $\varepsilon_n$ was arbitrary $B_\varepsilon(f(a), f(b))$ holds as required.  □

**Theorem 16.** *If $f \colon X \to Y$ where $X$ is a prelength space and $Y$ is a metric space and $x \colon \mathfrak{C}(X)$, then $\mathrm{map}(f)(x)$ is a regular function.*

**Proof.** Let $y \stackrel{def}{=} \mathrm{map}(f)(x)$. Consider arbitrary $\varepsilon_1$ and $\varepsilon_2$. Let $\delta_i \stackrel{def}{=} \mu_f(\varepsilon_i)$. Because $x$ is a regular function, $B_{\delta_1+\delta_2}(x(\delta_1), x(\delta_2))$ holds. By Lemma 15, $B_{\varepsilon_1+\varepsilon_2}(f(x(\delta_1)), f(x(\delta_2)))$ holds, which is the same as $B_{\varepsilon_1+\varepsilon_2}(y(\varepsilon_1), y(\varepsilon_2))$. Therefore, $y$ is a regular function.  □

**Theorem 17.** *The function $\mathrm{map}(f)$ is uniformly continuous with modulus $\mu_f$.*

**Proof.** Consider arbitrary $\varepsilon_0$ and $\varepsilon_1$. Let $\delta_i \stackrel{def}{=} \mu_f(\varepsilon_i)$. Suppose $B_{\delta_1}(x_1, x_2)$ holds. Thus $B_{\delta_0+\delta_1+\delta_2}(x_1(\delta_1), x_2(\delta_2))$ holds. By Lemma 15, $B_{\varepsilon_0+\varepsilon_1+\varepsilon_2}(f(x_1(\delta_1)), f(x_2(\delta_2)))$ holds. Therefore, $B_{\varepsilon_0}(\mathrm{map}(f)(x_1), \mathrm{map}(f)(x_2))$ holds.  □

The function bind : $(X \to \mathfrak{C}(Y)) \Rightarrow (\mathfrak{C}(X) \to \mathfrak{C}(Y))$ is defined in terms of join and map in the usual way.

$$\mathrm{bind}(f) \stackrel{def}{=} \mathrm{join} \circ \mathrm{map}(f).$$

## 4.1. Monad Laws.

To prove that completion forms a monad, one needs to verify that the seven monad laws [17] hold.

**Theorem 18.** $\mathrm{map}(\mathrm{id}) \asymp \mathrm{id}$, $\mathrm{map}(f \circ g) \asymp \mathrm{map}(f) \circ \mathrm{map}(g)$ *and* $\mathrm{map}(f) \circ \mathrm{unit} \asymp \mathrm{unit} \circ f$.

**Proof.** Trivial.  □

**Theorem 19.** $\mathrm{map}(f) \circ \mathrm{join} \asymp \mathrm{join} \circ \mathrm{map}(\mathrm{map}(f))$.

**Proof.** Let $x$, $\varepsilon_1$, and $\varepsilon_2$ be arbitrary. Let $\delta_1 \stackrel{def}{=} \mu_f(\varepsilon_1)$, $\delta_2 \stackrel{def}{=} \mu_f(\frac{\varepsilon_2}{2})$, $a \stackrel{def}{=} x(\frac{\delta_1}{2})(\frac{\delta_1}{2})$, and $b \stackrel{def}{=} x(\delta_2)(\delta_2)$. By expanding definitions, $\mathrm{map}(f)(\mathrm{join}(x))(\varepsilon_1) = f(a)$ and $\mathrm{join}(\mathrm{map}(\mathrm{map}(f))(x))(\varepsilon_2) = f(b)$. Because $x$ is a regular function, $B_{\frac{\delta_1}{2}+\delta_2}(x(\frac{\delta_1}{2}), x(\delta_2))$ and $B_{\delta_1+2\delta_2}(a, b)$ hold. By Lemma 15, $B_{\varepsilon_1+\varepsilon_2}(f(a), f(b))$ holds as required.  □

**Theorem 20.** $\mathrm{join} \circ \mathrm{unit} \asymp \mathrm{id}$.

**Proof.** Let $x$, $\varepsilon_1$, and $\varepsilon_2$, be arbitrary. $\mathrm{join}(\mathrm{unit}(x))(\varepsilon_1) = x(\frac{\varepsilon_1}{2})$ and $\mathrm{id}(x)(\varepsilon_2) = x(\varepsilon_2)$. Since $x$ is a regular function, $B_{\varepsilon_1 + \varepsilon_2}(x(\frac{\varepsilon_1}{2}), x(\varepsilon_2))$ holds as required. $\qquad\square$

**Theorem 21.** $\mathrm{join} \circ \mathrm{map}(\mathrm{unit}) \asymp \mathrm{id}$.

**Proof.** $\mathrm{join}(\mathrm{map}(\mathrm{unit})(x))(\varepsilon) = \mathrm{join}(\mathrm{unit}(x))(\varepsilon)$; therefore, this theorem follows from Theorem 20. $\qquad\square$

**Theorem 22.** $\mathrm{join} \circ \mathrm{map}(\mathrm{join}) \asymp \mathrm{join} \circ \mathrm{join}$.

**Proof.** By unfolding definitions, $\mathrm{join}(\mathrm{map}(\mathrm{join})(x))(\varepsilon_1) = x(\frac{\varepsilon_1}{2})(\frac{\varepsilon_1}{4})(\frac{\varepsilon_1}{4})$ and $\mathrm{join}(\mathrm{join}(x))(\varepsilon_2) = x(\frac{\varepsilon_2}{4})(\frac{\varepsilon_2}{4})(\frac{\varepsilon_2}{2})$. Because $x$ is a regular function, $B_{\frac{\varepsilon_1}{2} + \frac{\varepsilon_2}{4}}(x(\frac{\varepsilon_1}{2}), \quad x(\frac{\varepsilon_2}{4}))$, $B_{\frac{3\varepsilon_1}{4} + \frac{\varepsilon_2}{2}}(x(\frac{\varepsilon_1}{2})(\frac{\varepsilon_1}{4}), x(\frac{\varepsilon_2}{4})(\frac{\varepsilon_2}{4}))$, and $B_{\varepsilon_1 + \varepsilon_2}(x(\frac{\varepsilon_1}{2})(\frac{\varepsilon_1}{4})(\frac{\varepsilon_1}{4}), x(\frac{\varepsilon_2}{4})(\frac{\varepsilon_2}{4})(\frac{\varepsilon_2}{2}))$ hold. This last one is what is required. $\qquad\square$

### 4.2. Lifting Binary Functions.

The function map can be used to lift unary functions from $X \to Y$ to $\mathfrak{C}(X) \to \mathfrak{C}(Y)$; however, one also wants to lift multi-argument functions such as binary functions. One can lift curried binary functions on base spaces to curried functions on completed spaces, but one first needs to consider the space of uniformly continuous functions as a metric space. This is done using the supremum norm.

**Theorem 23.** *For any two metric spaces* $(X, \asymp, B^X)$ *and* $(Y, \asymp, B^Y)$, *The space of uniformly continuous functions from $X$ to $Y$, $(X \to Y, \asymp, B^{X \to Y})$, is a metric space where* $B_\varepsilon^{X \to Y}(f, g) \overset{\text{def}}{=} \forall a, B_\varepsilon^Y(f(a), g(a))$ *and with the equivalence relation* $f \asymp g \overset{\text{def}}{=} \forall a, f(a) \asymp g(a)$.

**Proof.** All the properties of $B^{X \to Y}$ follow directly from the same properties for $B^Y$. $\qquad\square$

Unfortunately, the space of uniformly continuous functions between two prelength spaces is not necessarily a prelength space. Consider the space of curves on the unit circle, $[0, 1] \to S^1$. Consider two curves connecting two points where one goes clockwise and the other goes counter-clockwise. If $[0, 1] \to S^1$ were a prelength space, then there would be a curve about halfway between these two curves; however, there is no such curve.

Fortunately, one can still lift curried functions using map, because the proof that the result of map is a regular function (see Theorem 16) requires only that the *domain* be a prelength space. With curried functions (say $X \to (X \to X)$), it is the *range* that might not be a prelength space. The domain is still the prelength space $X$, so Theorem 16 still applies.

One can show that the monad is a strong monad by showing that map is uniformly continuous.

**Lemma 24.** *Let $X$ be a prelength space and $Y$ a metric space. Then* map $: (X \to Y) \to (\mathfrak{C}(X) \to \mathfrak{C}(Y))$ *is uniformly continuous with modulus* id.

**Proof.** Let $\varepsilon$ be arbitrary. Suppose $f g : X \to Y$ are such that $B_\varepsilon(f, g)$. Let $\varepsilon_0$ and $x : \mathfrak{C}(X)$ be arbitrary. Let $\delta_0 \overset{\text{def}}{=} \min(\mu_f(\varepsilon_0), \mu_g(\varepsilon_0))$ and $a_0 \overset{\text{def}}{=} x(\delta_0)$. $B_{\varepsilon_0}(\mathrm{map}(f)(x), \mathrm{map}(f)(\mathrm{unit}(a_0)))$ and $B_{\varepsilon_0}(\mathrm{map}(g)(x), \mathrm{map}(g)(\mathrm{unit}(a_0)))$ both hold by uniform continuity and Theorem 13. By Theorem 18, $\mathrm{map}(f)(\mathrm{unit}(a_0)) \asymp \mathrm{unit}(f(a_0))$ and $\mathrm{map}(g)(\mathrm{unit}(a_0)) \asymp \mathrm{unit}(g(a_0))$. By Theorem 12, $B_\varepsilon(\mathrm{unit}(f(a_0)), \mathrm{unit}(g(a_0)))$ holds. So all together $B_{\varepsilon + 2\varepsilon_0}(\mathrm{map}(f)(x), \mathrm{map}(g)(x))$. Since $\varepsilon_0$ and $x$ were arbitrary, $B_\varepsilon(\mathrm{map}(f), \mathrm{map}(g))$ holds. $\qquad\square$

With uniformly continuous function spaces in hand, one can construct the function ap $: \mathfrak{C}(X \to Y) \to \mathfrak{C}(X) \to \mathfrak{C}(Y)$, which one can think of as a function that takes the limit of uniformly continuous functions.

$$\mathrm{ap}(f) \overset{\text{def}}{=} \lambda x. \lambda \varepsilon. \, \mathrm{map}(f(\tfrac{\varepsilon}{2}))(x)(\tfrac{\varepsilon}{2}).$$

**Theorem 25.** *For any* $f : \mathfrak{C}(X \to Y)$ *and* $x : \mathfrak{C}(X)$, $\mathrm{ap}(f)(x)$ *is a regular function.*

**Proof.** Let $\varepsilon_1$ and $\varepsilon_2$ be arbitrary. Let $f_i \overset{def}{=} f(\frac{\varepsilon_i}{2}) : X \to Y$. Because $f$ is regular, $B_{\frac{\varepsilon_1}{2} + \frac{\varepsilon_2}{2}}(f_1, f_2)$ holds. Let $y_i \overset{def}{=} \mathrm{map}(f_i)(x)$. By Lemma 24, $B_{\frac{\varepsilon_1}{2} + \frac{\varepsilon_2}{2}}(y_1, \quad y_2)$ holds. By Theorem 13, $B_{\frac{\varepsilon_1}{2}}(\mathrm{unit}(y_1(\frac{\varepsilon_1}{2})), \quad y_1)$ and $B_{\frac{\varepsilon_2}{2}}(y_2, \quad \mathrm{unit}(y_2(\frac{\varepsilon_2}{2})))$ hold. All together $B_{\varepsilon_1 + \varepsilon_2}(\mathrm{unit}(y_1(\frac{\varepsilon_1}{2})), \mathrm{unit}(y_2(\frac{\varepsilon_2}{2})))$ holds. By Theorem 12, $B_{\varepsilon_1 + \varepsilon_2}(y_1(\frac{\varepsilon_1}{2}), y_2(\frac{\varepsilon_2}{2}))$ as required. $\square$

**Lemma 26.** *For any* $f : \mathfrak{C}(X \to Y)$, $x : \mathfrak{C}(X)$, *and* $\varepsilon$, $B_\varepsilon(\mathrm{ap}(f)(x), \mathrm{map}(f(\varepsilon))(x))$ *holds.*

**Proof.** Let $\delta$ be arbitrary. Let $g \overset{def}{=} f(\frac{\delta}{2})$ and $h \overset{def}{=} f(\varepsilon)$. Because $f$ is regular, $B_{\frac{\delta}{2} + \varepsilon}(\mathrm{map}(g)(x), \mathrm{map}(h)(x))$ holds. By Theorem 13, $B_{\frac{\delta}{2}}(\mathrm{unit}(\mathrm{map}(g)(x)(\frac{\delta}{2})), \mathrm{map}(g)(x))$ and $B_\delta(\mathrm{map}(h)(x), \mathrm{unit}(\mathrm{map}(h)(x)(\delta)))$ hold. All together $B_{\varepsilon + 2\delta}(\mathrm{unit}(\mathrm{map}(g)(x)(\frac{\delta}{2})), \mathrm{unit}(\mathrm{map}(h)(x)(\delta)))$ holds, and by Theorem 12, $B_{\varepsilon + 2\delta}(\mathrm{map}(g)(x)(\frac{\delta}{2}), \mathrm{map}(h)(x)(\delta))$ holds. By Theorem 8, $B_\varepsilon(\mathrm{map}(g)(x), \mathrm{map}(h)(x))$ holds as required. $\square$

**Theorem 27.** *For any* $f : \mathfrak{C}(X \to Y)$, $\mathrm{ap}(f)$ *is uniformly continuous with modulus* $\lambda \varepsilon . \mu_{f(\frac{\varepsilon}{3})}(\frac{\varepsilon}{3})$.

**Proof.** Let $\varepsilon$ be arbitrary, and let $\delta \overset{def}{=} \mu_{f(\frac{\varepsilon}{3})}(\frac{\varepsilon}{3})$. Suppose $x_1$ and $x_2$ are such that $B_\delta(x_1, x_2)$ holds. Then $B_{\frac{\varepsilon}{3}}(\mathrm{map}(f(\frac{\varepsilon}{3}))(x_1), \mathrm{map}(f(\frac{\varepsilon}{3}))(x_2))$ holds. By Lemma 26, $B_{\frac{\varepsilon}{3}}(\mathrm{ap}(f)(x_i), \mathrm{map}(f(\frac{\varepsilon}{3}))(x_i))$ hold. Together this means that $B_\varepsilon(\mathrm{ap}(f)(x_1), \mathrm{ap}(f)(x_2))$ holds as required. $\square$

**Theorem 28.** *The function* $\mathrm{ap}$ *is uniformly continuous with modulus* $\mathrm{id}$.

**Proof.** Let $\varepsilon$ be arbitrary. Suppose $f_1$ and $f_2$ are such that $B_\varepsilon(f_1, f_2)$ holds. Let $x$ and $\delta$ be arbitrary. Because $f_1$ and $f_2$ are regular functions, $B_{\frac{\delta}{2} + \varepsilon}(f_1(\frac{\delta}{4}), f_2(\frac{\delta}{4}))$ holds. Therefore, $B_{\frac{\delta}{2} + \varepsilon}(\mathrm{map}(f_1(\frac{\delta}{4}))(x), \mathrm{map}(f_2(\frac{\delta}{4}))(x))$ holds. By Lemma 26, $B_{\frac{\delta}{4}}(\mathrm{ap}(f_i)(x), \mathrm{map}(f_i(\frac{\delta}{4}))(x))$ holds. Together this yields $B_{\delta + \varepsilon}(\mathrm{ap}(f_1)(x), \mathrm{ap}(f_2)(x))$, and since $\delta$ and $x$ were arbitrary, $B_\varepsilon(\mathrm{ap}(f_1), \mathrm{ap}(f_2))$ holds as required. $\square$

By using $\mathrm{ap}$, one can define a $\mathrm{map2}$ function that lifts a curried function $f : X \to Y \to Z$.

$$\mathrm{map2}(f) \overset{def}{=} \mathrm{ap} \circ \mathrm{map}(f) : \mathfrak{C}(X) \to \mathfrak{C}(Y) \to \mathfrak{C}(Z).$$

Other map functions can be made to lift $n$-ary functions by repeated use of $\mathrm{ap}$.

## 5. Functions of Real Numbers

Using this completion monad, it is now very easy to define the real numbers.

$$\mathbb{R} \overset{def}{=} \mathfrak{C}(\mathbb{Q})$$

Functions operating on the real numbers can be created by using $\mathrm{map}$ or $\mathrm{bind}$ to lift uniformly continuous functions operating on the rationals. This is the advantage of using the monad functions. It is often easier to define a function on the rationals than on the reals because it is possible to decide for two rational numbers $a$ and $b$ which of $a < b$, $a = b$, or $a > b$ hold. The real numbers do not have this property.

**Theorem 29.**

- For all $a : \mathbb{Q}$, the functions $\lambda b. - b$ and $\lambda b. a + b$ are both uniformly continuous from $\mathbb{Q}$ to $\mathbb{Q}$ with moduli $\mathrm{id}$.

- For all $a : \mathbb{Q}$, $\lambda b. \max(a, b)$ is uniformly continuous from $\mathbb{Q}$ to $[a, \infty[$ with modulus $\mathrm{id}$.

- For all $a : \mathbb{Q}$, $\lambda b. \min(a, b)$ is uniformly continuous from $\mathbb{Q}$ to $]-\infty, a]$ with modulus $\mathrm{id}$.

- The function $\lambda b. |b|$ is uniformly continuous from $\mathbb{Q}$ to $[0, \infty[$ with modulus $\mathrm{id}$.

- *For all $a: \mathbb{Q}$, $\lambda b. \varepsilon\, a\, b$ is uniformly continuous from $\mathbb{Q}$ to $\mathbb{Q}$ with modulus $\lambda \varepsilon. \varepsilon\; |a|^{-1}$ (or any modulus in the case when $a = 0$).*

All these unary functions can be lifted, by using map, to functions on $\mathbb{R}$. For the moment one of the arguments to the binary functions $+\,$, $\cdot\,$, $\min\,$, and $\max$ must be rational. For example, for any $a: \mathbb{Q}$ one can use map to create the unary function $\lambda b. a + b: \mathbb{R} \to \mathbb{R}$, the function that translates by a rational number, but one does not yet have the binary function $\lambda a\, b. a\; +\; b:$ $\mathbb{R} \to \mathbb{R} \to \mathbb{R}$.

At this point, one can see why regular functions work better than regular sequences. For multiplication the modulus of continuity is $\lambda \varepsilon. \varepsilon\; |a|^{-1}$. If one used regular sequences, then one would need to find a point suitably far in the sequence. The result would be an over-approximation of the computation one needs. For a deep computation, these over-approximations would add up and can cause a lot of unnecessary work to be done. By using regular functions, one can request exactly the amount of accuracy required.

### 5.1. Reciprocal.

Not all functions that one wishes to represent are uniformly continuous. One example of a non-uniformly continuous function is the reciprocal function, $\lambda x. x^{-1}$. Bishop and Bridges [1] define a continuous function as a function that is uniformly continuous on every closed interval in its domain. The general idea is for each $x$ to find a domain containing $x$ that the function is uniformly continuous on, and then evaluate this function at $x$.

**Theorem 30.** *For every $a: \mathbb{Q}^+$, $\lambda x. x^{-1}$ is uniformly continuous on $[a, \infty[$ and $]\infty, -a]$, and both functions have a modulus of continuity of $\lambda \varepsilon. \varepsilon\, a^2$.*

One can now define $x^{-1}$ for any $x: \mathbb{R}$ such that $x < 0$ or $0 < x$, where $0 < x \overset{def}{=} \exists \varepsilon. \varepsilon < x(\varepsilon)$ and similarly for $x < 0$. Suppose $0 < x$. Then there exists some $a: \mathbb{Q}$ such that $0 < a \le x$. One can inject $x$ into the domain $[a, \infty[$, and then invert it.

$$x^{-1} \overset{def}{=} \mathrm{map}(\lambda b: [a, \infty[.b^{-1})(\max(a, x)).$$

Similarly if $x < 0$, then there exists some $a: \mathbb{Q}$ such that $x \le a < 0$. In this case, define the inverse as

$$x^{-1} \overset{def}{=} \mathrm{map}(\lambda b: ]-\infty, a].b^{-1})(\min(a, x)).$$

It is important to note that the result is independent of the choice of domain that $x$ is injected into. If $0 < a \le x$ and $0 < a' \le x$, then whether $x$ is injection into $[a, \infty[$ or $x$ is injected into $[a', \infty[$, the result of $x^{-1}$ is equivalent.

Notice that the larger the value of $a$, the larger the modulus of continuity is. The larger the modulus of continuity is, the more efficient the lifted function is because less work is done in approximating the input. Therefore, it is helpful to find as large an $a$ as is reasonable when computing the reciprocal.

### 5.2. Binary Functions of Real Numbers.

One can use map2 to lift addition to operate on $\mathbb{R}$.

**Theorem 31.** *The function $\lambda a. \lambda b. a\; +\; b$ is uniformly continuous from $\mathbb{Q}$ to $\mathbb{Q} \to \mathbb{Q}$ with modulus id.*

The definition of addition for $x\, y: \mathbb{R}$ is

$$x + y \overset{def}{=} \mathrm{map2}(\lambda a. \lambda b. a + b)(x)(y).$$

The definition of multiplication on $\mathbb{R}$ is a little more complicated than for addition.

**Theorem 32.** *For any $c: \mathbb{Q}^+$, the function $\lambda a. \lambda b. a\, b$ is uniformly continuous from $\mathbb{Q}$ to $[-c, c] \to \mathbb{Q}$ with modulus $\lambda \varepsilon. \varepsilon\, c^{-1}$.*

The definition of multiplication for $x\,y\colon\mathbb{R}$ is

$$x\,y \overset{\text{def}}{=} \text{map2}(\lambda a.\lambda b.\,a\,b)(x)(\max{(-c,\min{(c,y)})}) \qquad \text{where } c \overset{\text{def}}{=} |y(1)| + 1.$$

Because multiplication is not uniformly continuous, some closed interval containing $y$ must be found. This requirement is reflected in the fact that $c$ is used in the definition of the modulus of continuity.

This definition of multiplication is asymmetric because only $y$ is bound by a closed interval. The usual definition of multiplication in exact real arithmetic implementations finds closed intervals containing both $x$ and $y$. The definition here is the natural definition of multiplication one gets when one uses monad operations. It saves one from the unnecessary work of finding a closed interval for $x$.

### 5.3. Power Series.

Other elementary functions can be defined in terms of power series. The limit of a convergent sequence can be found if the modulus of convergence is known. It is particularly easy to compute the limit of an alternating decreasing series.

**Theorem 33.** *Suppose* $\lim_{n \to \infty} a_n = 0$, $a_n$ *is alternating, and* $\forall i,\ |a_{i+1}| < |a_i|$. *Then* $\sum_{i=0}^{\infty} a_i = x$ *where* $x(\varepsilon) \overset{\text{def}}{=} \sum_{\varepsilon < a_i} a_i$.

The power series of several elementary functions are alternating and decreasing on small rational inputs because the terms are bounded by a geometric series.

- For $-1 < a \le 0$, $\exp^{\mathbb{Q}}(a) \overset{\text{def}}{=} \sum_{i=0}^{\infty} \frac{a^i}{i!}$.

- For $-1 < a < 1$, $\sin^{\mathbb{Q}}(a) \overset{\text{def}}{=} \sum_{i=0}^{\infty} (-1)^i \frac{a^{2i+1}}{(2i+1)!}$.

- For $1 \le a < 2$, $\ln^{\mathbb{Q}}(a) \overset{\text{def}}{=} \sum_{i=0}^{\infty} (-1)^i \frac{(a-1)^{i+1}}{i+1}$.

- For $-1 < a < 1$, $\arctan^{\mathbb{Q}}(a) \overset{\text{def}}{=} \sum_{i=0}^{\infty} (-1)^i \frac{a^{2i+1}}{2i+1}$.

These four functions can be extended to all rationals in their domains by repeated applications of the following formulas to reduce input into the above domains.

$$\exp(a) \;\asymp\; \exp^2(\frac{a}{2})$$
$$\exp(a) \;\asymp\; \exp^{-1}(-a)$$
$$\sin(a) \;\asymp\; 3\sin(\frac{a}{3}) - 4\sin^3(\frac{a}{3})$$
$$\ln(a) \;\asymp\; \ln(\frac{a}{2}) + \ln(2)$$
$$\ln(a) \;\asymp\; \ln(\frac{3}{4}a) + \ln(\frac{4}{3})$$
$$\ln(a) \;\asymp\; -\ln(a^{-1})$$
$$\arctan(a) \;\asymp\; \arctan(-a)$$
$$\arctan(a) \;\asymp\; \frac{\pi}{2} - \arctan(a^{-1}) \qquad (\text{for } 0 < a)$$
$$\arctan(a) \;\asymp\; \frac{\pi}{4} + \arctan(\frac{a-1}{a+1}) \qquad (\text{for } 0 < a)$$

The number $\pi$ can be defined in terms of arctan in several ways. Below is a definition [7] that uses particularly small inputs to arctan.

$$\pi \overset{\text{def}}{=} \left(48\arctan(\frac{1}{38}) + 80\arctan(\frac{1}{57})\right) + \left(28\arctan(\frac{1}{239}) + 96\arctan(\frac{1}{268})\right)$$

Also cos can be easily defined in terms of sin.

$$\cos(a) \overset{\text{def}}{=} 1 - 2\sin^2(\frac{a}{2})$$

**Theorem 34.** *The functions* sin, cos, *and* arctan *are uniformly continuous on* $\mathbb{Q}$ *with moduli* id.

**Theorem 35.** *For every $a : \mathbb{Z}$, exp is uniformly continuous on $]-\infty, a]$ with modulus $\lambda \varepsilon. \varepsilon\, 2^{-a}$ if $a \leq 0$ and modulus $\lambda \varepsilon. \varepsilon\, 3^{-a}$ if $0 \leq a$.*

**Theorem 36.** *For every $a : \mathbb{Q}^{+}$, ln is uniformly continuous on $[a, \infty[$ with modulus $\lambda \varepsilon. \varepsilon\, a$.*

These functions can all be lifted, using bind, to functions on $\mathbb{R}$ or whatever the completions of their domains are. The non-uniformly continuous functions exp and ln are defined the same way the reciprocal was defined; for each input $x$ some suitable domain is found, and $x$ is injected into that domain. Then the lifted function is evaluated.

Other elementary functions can be defined in terms of the four elementary functions given above.

$$\tan(x) \stackrel{def}{=} \frac{\sin(x)}{\cos(x)}$$

$$\sinh(x) \stackrel{def}{=} \frac{\exp(x) - \exp(-x)}{2}$$

$$\cosh(x) \stackrel{def}{=} \frac{\exp(x) + \exp(-x)}{2}$$

$$\tanh(x) \stackrel{def}{=} \frac{\sinh(x)}{\cosh(x)}$$

$$x^{y} \stackrel{def}{=} \exp(y \ln(x))$$

$$\arcsin(x) \stackrel{def}{=} \arctan(\frac{x}{\sqrt{1 - x^2}})$$

$$\arccos(x) \stackrel{def}{=} \frac{\pi}{2} - \arcsin(x)$$

$$\operatorname{arcsinh}(x) \stackrel{def}{=} \ln(x + \sqrt{1 + x^2})$$

$$\operatorname{arccosh}(x) \stackrel{def}{=} \ln(x + \sqrt{x^2 - 1})$$

$$\operatorname{arctanh}(x) \stackrel{def}{=} \frac{1}{2} \ln(\frac{1 + x}{1 - x})$$

## 6. Improving Efficiency

### 6.1. Compression.

Calculations on rational numbers become more expensive as the size of the numerator and denominator become larger. One can keep rational numbers in reduced form by dividing through by the gcd of the numerator and denominator after each operation; however, for real numbers one can do better. Let us say that $\frac{n_1}{d_1}$ is simpler than $\frac{n_2}{d_2}$ if $|n_1| + |d_1| < |n_2| + |d_2|$. In any closed rational interval, there is a unique rational number that is simplest [16]. This idea can be used to *compress* real numbers; that is to try to make the rational outputs of regular functions simple while keeping the real number it represents the same.

$$\operatorname{approx}(a, \varepsilon) \stackrel{def}{=} c \qquad \text{where } c \text{ is the simplest rational number in } [a - \varepsilon, a + \varepsilon].$$

$$\operatorname{compress}(x) \stackrel{def}{=} \lambda \varepsilon. \operatorname{approx}(x(\frac{\varepsilon}{2}), \frac{\varepsilon}{2}).$$

**Theorem 37.** *For any $x : \mathbb{R}$, $\operatorname{compress}(x) \asymp x$.*

**Proof.** Let $\varepsilon$ be arbitrary. $B_{\frac{\varepsilon}{2}}(\operatorname{compress}(x)(\varepsilon), x(\frac{\varepsilon}{2}))$ holds by the definition of approx. Because $x$ is a regular function, $B_{\frac{3\varepsilon}{2}}(x(\frac{\varepsilon}{2}), x(\varepsilon))$ holds. Therefore, $B_{\varepsilon}(\operatorname{compress}(x)(\varepsilon), x(\varepsilon))$ holds as required. $\qquad\square$

Generous use of compress can greatly increase the efficiency of the computation. In my implementation every function operating on $\mathbb{R}$ passes its input through compress first.

Notice that if one uses the alternative defintion of regular function given in section 3.4, then this definition of compress does not work. I do not see a way of defining a version of compress for the alternative definition of regular function, and without compress, the rational numbers used in computation become too large for practical evaluation.

### 6.2. Newton Approximation.

Although one can define $\sqrt{x}$ as $\exp(\frac{1}{2}\ln(x))$, the square root function is common enough to warrant its own, more efficient, definition. If $a\colon\mathbb{Q}$ is such that $1 \le a < 4$, one can compute $\sqrt{a}$ by taking the limit of Newton approximations. Let the first approximation be $b_0 \overset{\text{def}}{=} \frac{a+1}{2}$, and let successive approximations $b_{i+1} \overset{\text{def}}{=} \mathrm{approx}_{\mathbb{Q}}(\frac{a+b_i^2}{2\,b_i}, 2^{-2^{i+1}-1})$. The first approximation has an error of at most $\frac{1}{2}$. Each successive approximation squares the error.

**Theorem 38.** *For any* $a\colon\mathbb{Q}$, *such that* $1 \le a < 4$, *let* $x(\varepsilon) \overset{\text{def}}{=} b_n$, *where* $n$ *is the first natural number such that* $2^{-2^n} \le \varepsilon$. *Then* $x$ *is a regular function.*

**Proof.** It suffices to show that for all $n$, $\sqrt{a} - 2^{-2^n} \le b_n \le \sqrt{a} + 2^{-2^n}$ and $1 \le b_n$.

Suppose $n = 0$. For the bound below note that $0 \le \frac{(a-1)^2}{4} = b_0^2 - a$, so $\sqrt{a} - 2^{-1} \le \sqrt{a} \le b_0$. For the bound above note that $(b_0 - 2^{-1})^2 = \frac{a^2}{4} \le a$, so $b_0 \le \sqrt{a} + 2^{-1}$. Finally note that $1 \le a$, so $1 \le b_0$.

Suppose $n = m + 1$. Let $\delta \overset{\text{def}}{=} 2^{-2^m}$ and $c \overset{\text{def}}{=} \frac{a+b_m^2}{2\,b_m}$. Assume by induction that $\sqrt{a} - \delta \le b_m \le \sqrt{a} + \delta$ and $1 \le b_m$. This means that $|b_m - \sqrt{a}| \le \delta$. This implies that $c - \sqrt{a} \le \frac{\delta^2}{2\,b_m} \le \frac{\delta^2}{2}$ because $1 \le b_m$. Also notice that $0 \le \frac{(b_m + \sqrt{a})^2}{2\,b_m} = c - \sqrt{a}$. Together one sees that $\sqrt{a} \le c \le \sqrt{a} + \frac{\delta^2}{2}$. Because $b_{m+1}$ is within $\frac{\delta^2}{2}$ of $c$, one concludes that $\sqrt{a} - \delta^2 \le b_{m+1} \le \sqrt{a} + \delta^2$. All that remains is to prove that $1 \le b_{m+1}$. $1 \le c$ because $1 \le a$. The only rational simpler than 1 is 0. Therefore, the only way $b_{m+1}$ could be less than 1 is if it were 0. That would require that $1 \le c \le \frac{\delta^2}{2}$; however $\frac{\delta^2}{2}$ is less than 1. This makes it is impossible for $c$ to be 0. $\qquad\square$

For $0 < a$ that falls outside the interval $[1, 4[$, one can find some $m\colon\mathbb{Z}$ such that $1 \le 4^m\,a < 4$. Therefore, one can define the square root function on the rationals as

$$\sqrt{0} \overset{\text{def}}{=} 0$$
$$\text{for } 1 \le a < 4, \ \sqrt{a} \overset{\text{def}}{=} x \quad \text{(where } x \text{ is defined as in Theorem 38)}$$
$$\text{otherwise, } \sqrt{a} \overset{\text{def}}{=} \frac{\sqrt{4^m\,a}}{2^m} \quad \text{(for some } m\colon\mathbb{Z} \text{ such that } 1 \le 4^m\,a < 4).$$

**Theorem 39.** *The function* $\lambda a.\ \sqrt{a}$ *is uniformly continuous from* $\mathbb{Q}^{0+}$ *to* $\mathbb{R}$, *with modulus* $\lambda\varepsilon.\varepsilon^2$.

By using bind, the square root function can be lifted to the domain $\mathfrak{C}(\mathbb{Q}^{0+})$.

### 6.3. More Efficient Power Series.

A power series converges faster for values closer to 0. For sin and exp, the two equations mentioned before

$$\exp(a) \asymp \exp^2(\frac{a}{2})$$
$$\sin(a) \asymp 3\sin(\frac{a}{3}) - 4\sin^3(\frac{a}{3})$$

can be repeatedly applied to shrink the input arbitrarily close to 0. There is, of course, a trade-off between evaluating a polynomial each time this equation is used and having the power series converge faster. The optimal trade-off depends on the input and the implementation of the library. For simplicity, in my implementation I shrink the input until it is less than $2^{-50}$, which is a value that seems to work well in my implementation for some example problems [13].

Unfortunately such nice reductions do not seem to exist for ln or arctan.

### 6.4. More Efficient Periodic Functions.

The functions sin and cos are both periodic with period $2\pi$. With a reasonably fast implementation of $\pi$, it is possible to subtract out a multiple of $2\pi$ to reduce the size of the argument to sin or cos.

$$\forall n\colon\mathbb{Z}, \sin(x) \asymp \sin(x - n\,2\pi)$$
$$\forall n\colon\mathbb{Z}, \cos(x) \asymp \cos(x - n\,2\pi)$$

The best value of $n$ would be $\lfloor \frac{x}{2\pi} \rceil$. The function mapping $x$ to a nearest integer $\lfloor x \rceil$ is not computable for real numbers, but it is computable for the rationals—it does not matter which integer $\frac{1}{2}$ maps to. The approximation $\frac{x}{2\pi}(\frac{1}{2})$ is within $\frac{1}{2}$ of $\frac{x}{2\pi}$, so $\lfloor \frac{x}{2\pi}(\frac{1}{2}) \rceil$ is an integer within 1 of $\lfloor \frac{x}{2\pi} \rceil$, which is good enough to use for $n$.

In practice one can be even more clever with the symmetries of sin and cos to reduce the range a little bit further.

**6.5. Summing Lists.**

It is not uncommon to want to sum a list of real numbers. The use of ap in the definition of addition means that to compute $x + y$ within $\varepsilon$, one must compute both $x$ and $y$ within $\frac{\varepsilon}{2}$. If one associates the sum of a list of numbers all to one side $a_0 + (a_1 + (a_2 + ... + a_n))$, the first term $a_0$ is only approximated within $\frac{\varepsilon}{2}$, while the last term $a_n$ is approximated within $\frac{\varepsilon}{2^n}$.

If one is summing a list with exactly a power of 2 number of terms, as in the definition of $\pi$, then one can associate the sum into a balanced binary tree. In this case all terms will be approximated to the same degree; however, to handle more general cases it is useful to directly define the sum of $n$ real numbers.

$$\sum_{i=0}^{n} x_i \stackrel{def}{=} \lambda \varepsilon. \sum_{i=0}^{n} x_i(\frac{\varepsilon}{n})$$

**Theorem 40.** *The term $\sum_{i=0}^{n} x_i$ is a regular function.*

## 7. Implementation

The one goal of this work is to produce an actual implementation of the real numbers that is reasonably efficient and, more importantly, reliable. I have written a prototype Haskell module that implements the functions described in this paper. This code, entitled *Few Digits* [14], competed in the *"Many Digits" Friendly Competition* [13] occurring on October 3-4, 2005. The library consists of less than 440 lines of Haskell code and is included in Appendix A. Although the code did not do particularly well in the contest, finishing eighth out of nine, it still can compute thousands of decimal digits within minutes or seconds for many of the problems.

A lot of time is spent computing approx$(a, \varepsilon)$. One suggestion I have received is to instead quickly find some good approximation instead of spending so much time finding the perfect approximation. Another place much time is spent is evaluating the polynomial approximation of the elementary functions. It only takes a polynomial of a few degrees to cause the exact evaluation on large (in terms of simplicity) rational numbers to become very slow. Such precise computation over the rationals is unnecessary because the resulting rational number will just be feed to compress and a simpler rational number found. One possible solution may be to define the reals to be the completion of some dense set other than $\mathbb{Q}$ (such as the dyadic rationals $\mathbb{D} \stackrel{def}{=} \{\frac{a}{2^n} | a\, n : \mathbb{Z}\}$). However, such a set is not closed under division, so evaluating Taylor series becomes more difficult.

The code is more or less faithful to the theory described in this paper. In the Haskell code `Complete` is the type constructor of the monad, unit is called `const`, join is called `completeComplete`, map is called `evalUniformCts`, and bind is called `evalUniformCts2`. Unfortunately, because only uniformly continuous functions can be mapped, this monad does not fit into Haskell's `Monad` class.

The most significant difference between the description here and the implementation is that the power series for exp is computed for small positive inputs. Even though the series is not alternating, the series is still nice enough to get a good handle on its modulus of convergence. Another minor difference is that map2, called `evalBinaryUniformCts`, is implemented directly rather than implemented via ap.

There are a couple of other Haskell specific issues to note. In order to make the reals an instance of the `Floating` class, it is necessary to make a new data type. This new data type is called `CReal`, which means *constructive reals* (or *computable reals* or *Cauchy reals*). Since a new data type must be created anyway, I also added a memoized integer approximation to the data structure.

## 7.1. Comparisons With Other Implementations.

Boehm *et al.* [2] were the first to implement an exact real arithmetic library. They implemented both a functional representation and a representation as lazy lists. The real number representation in this paper is a functional representation. Lester has very short Haskell module (a little over 200 lines) called Era (exact real arithmetic) that has a real number data type very similar to Boehm *et al.*'s implementation. Computer proofs of some of the algorithms have been verified by PVS [9][8]. Because of the similarity between the two programs, a quick comparison is in order (unfortunately Era did not compete in the *"Many Digits" Friendly Competition*). Here is one test that illustrates the difference between the two implementations. Consider computing an approximation of $\cos(\frac{F_4}{F_5})$ and $\cos(\frac{F_{2394}}{F_{2395}})$ within $10^{-1000}$ where $F_n$ is the $n^{\text{th}}$ Fibonacci number $(\frac{F_4}{F_5} = \frac{3}{5})$. The two ratios are both close to the golden ratio, so they are close to each other. Table 1 shows the time each system takes to do the real number portion of the computation (the Fibonacci ratios were precomputed).

|             | $\cos(\frac{F_4}{F_5})$ | $\cos(\frac{F_{2394}}{F_{2395}})$ |
|-------------|-------------------------|-----------------------------------|
| Few Digits  | 0.93 s                  | 37.53 s                           |
| Era         | 2.73 s                  | 2.71 s                            |

**Table 1.** Timings for the computation of two Haskell real number implementations.

Because Era uses dyadic rationals in its representation, the two ratios are equally bad for it. Era takes about the same amount of time for both problems. The power series implementation in Few Digits does polynomial computation with rational numbers before compressing. The more terms and the larger the rational numbers, the longer the computation takes.

Another comparable implementation is Muñoz and Lester's [12] PVS implementation. It uses interval arithmetic to perform real number computation. Dealing with non-monotonic functions is difficult with interval arithmetic. Because of this, their implementation cannot approximate some expressions, such as $\sin(\frac{\pi}{2})$, to arbitrary accuracy. The approach in this paper has no problems dealing with non-monotonic functions.

## 8. Future Work

In order to produce a highly reliable implementation of the real numbers, the next step is to write the proofs in a computer proof assistant such as Coq [10] and verify the correctness of the algorithms. Coq is a proof assistant based on the calculus of inductive constructions. It can be seen as a dependently typed functional programming language with inductive and coinductive data types. I believe that the theorems presented in this paper are sufficiently straight-forward to be easily proved in a computer proof assistant.

## 9. Conclusion

The completion monad offers a simple, clear, and convenient way of constructively creating the completion of any metric space, and defining continuous functions on prelength spaces. I have created one possible implementation of the real numbers using this style, and demonstrate that this approach can perform practical computations. This approach avoids some difficulties that other computer verified implementations have. This completion monad is flexible enough to allow for alternative implementations of real numbers, such as the completion of the dyadic rationals, and alternative implementations of elementary functions. This same monad could be used to create other complete spaces, such as the complex numbers, the unit sphere, the projective plane, etc.

## 10. Acknowledgements

I would like to thank Milad Niqui for convincing me to write a Haskell implementation for his contest. I would like to thank Freek Wiedijk for advising me on using Newton's method to compute square roots. Most importantly, I would like to thank my advisor Bas Spitters for his advice when I was developing these ideas.

This document has been produced using $\TeX_{\text{MACS}}$.

# Appendix A. Haskell Source Code for *Few Digits*

```
module Complete where

type Gauge = Rational -- Intended to be strictly positive

-- Intended that d (f x) (f y) <= x + y
type Complete a = Gauge -> a

completeComplete :: (Complete (Complete a)) -> Complete a
completeComplete f eps = (f (eps/2)) (eps/2)

-- A uniformly continuous function on some subset of a to b
-- Hopefully the name of the function gives an indication of
-- the domain.
data UniformCts a b = UniformCts
        {modulus :: (Gauge -> Gauge),
         forgetUniformCts :: (a -> b)}

evalUniformCts :: UniformCts a b -> Complete a -> Complete b
evalUniformCts (UniformCts mu f) x eps = f (x (mu eps))

evalUniformCts2 :: UniformCts a (Complete b) -> Complete a -> Complete b
evalUniformCts2 f x = completeComplete $ evalUniformCts f x

completeUniformCtsRange :: Complete (UniformCts a b) ->
                           UniformCts a (Complete b)
completeUniformCtsRange f = UniformCts mu g
 where
  mu eps = modulus (f (eps/2)) (eps/2)
  g x eps = forgetUniformCts (f eps) x

evalBinaryUniformCts :: UniformCts a (UniformCts b c) ->
                        Complete a -> Complete b -> Complete c
{- This was the original implementation.
evalBinaryUniformCts f x =
  evalUniformCts2 $ completeUniformCtsRange $ evalUniformCts f x
-}

{- This implementation seems better -}
evalBinaryUniformCts f x y eps = (evalUniformCts approxf y) (eps/2)
 where
  approxf = (evalUniformCts f x) (eps/2)

module Base where

import Ratio

type Base = Rational

approxBase :: Base -> Rational -> Base
approxBase = approxRational

powers x = zipWith (%) (genpowers n) (genpowers d)
   where
    n = numerator x
    d = denominator x
    genpowers x = 1:(map (x*) (genpowers x))

module CReal where

import Base
import Complete
import Control.Exception

radius :: Base
radius = 2^^(-51)

data CReal = CReal {approx :: Complete Base,
                    intApprox :: Integer}

bound x = fromInteger (1+(abs (intApprox x)))

makeCReal :: Complete Base -> CReal
makeCReal x = CReal x' n
 where
  n = round (x (1/2))
  x' eps | eps >= 1 = fromInteger n
         | otherwise = x eps

{- produces a regular function whose resulting approximations are
   small in memory size -}
compress :: Complete Base -> Complete Base
compress x eps = approxBase (x (eps/2)) (eps/2)

squish :: CReal -> Complete Base
squish = compress . approx

instance Show CReal where
 show x = error "Cannot show a CReal"
```

```
{- show x = show $ map (\n -> squish x ((1/2)^n)) [0..] -}

realBase :: Base -> CReal
realBase x = makeCReal (const x)

approxRange :: CReal -> Gauge -> (Base, Base)
approxRange x eps = (r-eps, r+eps)
 where
  r = approx x eps

{- proveNonZeroFrom will not terminate if the input is 0 -}
{- Finds some y st 0 < (abs y) <= (abs x) -}
proveNonZeroFrom :: Gauge -> CReal -> Base
proveNonZeroFrom g r | high < 0 = high
                     | 0 < low = low
                     | otherwise = proveNonZeroFrom (g/2) r
 where
  (low, high) = approxRange r g

proveNonZero = proveNonZeroFrom 1

makeCRealFun :: (UniformCts Base Base) -> CReal -> CReal
makeCRealFun f x = makeCReal $ evalUniformCts f (squish x)

makeCRealFun2 :: (UniformCts Base (Complete Base)) -> CReal -> CReal
makeCRealFun2 f x = makeCReal $ evalUniformCts2 f (squish x)

makeCRealBinFun :: (UniformCts Base (UniformCts Base Base)) ->
                   CReal -> CReal -> CReal
makeCRealBinFun f x y =
  makeCReal $ evalBinaryUniformCts f (squish x) (squish y)

negateCts = UniformCts id negate

realNegate :: CReal -> CReal
realNegate = makeCRealFun negateCts

plusBaseCts a = UniformCts id (a+)
realTranslate a = makeCRealFun (plusBaseCts a)

plusCts :: UniformCts Base (UniformCts Base Base)
plusCts = UniformCts id plusBaseCts

realPlus :: CReal -> CReal -> CReal
realPlus = makeCRealBinFun plusCts

instance Eq CReal where
 a==b = 0==proveNonZero (realPlus a (realNegate b))

multBaseCts 0 = UniformCts (const 1) (const 0)
multBaseCts a = UniformCts mu (a*)
 where
  mu eps = eps/(abs a)

realScale :: Base -> CReal -> CReal
realScale 0 = \_ -> realBase 0
realScale a = makeCRealFun (multBaseCts a)

{- \x -> (\y -> (x*y)) is uniformly continuous on the domain (abs y) <= maxy -}
multUniformCts :: Base ->
                  UniformCts Base (UniformCts Base Base)
multUniformCts maxy = UniformCts mu multBaseCts
 where
  mu eps = assert (maxy>0) (eps/maxy)

{- We need to bound the value of x or y.  I think it is better to bound
   x so I actually compute y*x -}
realMult :: CReal -> CReal -> CReal
realMult x y = makeCRealBinFun (multUniformCts (bound x)) y x

absCts = UniformCts id abs

realAbs :: CReal -> CReal
realAbs = makeCRealFun absCts

instance Num CReal where
 (+) = realPlus
 (*) = realMult
 negate = realNegate
 abs = realAbs
 signum x = realScale (signum (proveNonZero x)) (realBase 1)
 fromInteger = realBase . fromInteger

{- domain is (-inf, nonZero) if nonZero < 0
   domain is [nonZero, inf) if nonZero > 0 -}
recipUniformCts :: Base -> UniformCts Base Base
recipUniformCts nonZero = UniformCts mu f
 where
  f a | 0 <= nonZero = recip (max nonZero a)
      | otherwise = recip (min a nonZero)
  mu eps = eps*(nonZero^2)

realRecipWitness :: CReal -> Base -> CReal
```

```
realRecipWitness x nonZero = makeCRealFun (recipUniformCts nonZero) x

realRecip :: CReal -> CReal
realRecip x = realRecipWitness x (proveNonZero x)

instance Fractional CReal where
 recip = realRecip
 fromRational = realBase . fromRational

intPowerCts _ 0 = UniformCts (const 1) (const 1)
intPowerCts maxx n = UniformCts mu (^n)
 where
  mu eps = assert (maxx > 0) $ eps/((fromIntegral n)*(maxx^(n-1)))

realPowerInt x n = makeCRealFun (intPowerCts (bound x) n) x

type Polynomial a = [a]

evalPolynomial :: (Num a) => Polynomial a -> a -> a
evalPolynomial [] x = 0
evalPolynomial (a:as) x = a + x*(evalPolynomial as x)

diffPolynomial :: (Num a) => Polynomial a -> Polynomial a
diffPolynomial p = zipWith (*) (tail p) (map fromInteger [1..])

polynomialUniformCts :: Base ->
                        Polynomial Base -> UniformCts Base Base
polynomialUniformCts _ [] = UniformCts (const 1) (const 0)
polynomialUniformCts maxx p |maxSlope==0 = UniformCts (const 1) (const (head p))
                            |otherwise = UniformCts mu (evalPolynomial p)
 where
  maxSlope = evalPolynomial (map abs (diffPolynomial p)) (max 1 maxx)
  mu eps = assert (maxSlope > 0) $ eps/maxSlope

realBasePolynomial :: Polynomial Base -> CReal -> CReal
realBasePolynomial p x =
 makeCRealFun (polynomialUniformCts (bound x) p) x

factorials = fact 1 1
 where
  fact i j = i:(fact (i*j) (j+1))

interleave [] _ = []
interleave (x:xs) l = x:(interleave l xs)

taylorApprox m p x eps =
  sum (zipWith (*) p (takeWhile highError preTerms))
 where
  preTerms = zipWith (/) (powers x) factorials
  highError t = m*(abs t) >= eps

{- only valid for x <= ln(2).  Works best for |x| <= 1/2 -}
rationalSmallExp :: Base -> CReal
rationalSmallExp x = assert ((abs x)<=(1/2)) $
  makeCReal $ expTaylorApprox
 where
  m | x <= 0 = 1
    | otherwise = 2
  expTaylorApprox eps =
    sum terms
   where
    terms = takeWhile highError $ zipWith (/) (powers x) factorials
    highError t = m*(abs t) >= eps

rationalExp :: Base -> Base -> CReal
rationalExp tol x | (abs x) <= tol = rationalSmallExp x
                  | otherwise = realPowerInt (rationalExp tol (x/2)) 2

expUniformCts :: Integer -> UniformCts Base (Complete Base)
expUniformCts upperBound = UniformCts mu (approx . rationalExp radius)
 where
  mu eps | upperBound <= 0 = eps*(2^(-upperBound))
         | otherwise = eps/(3^upperBound)

realExp :: CReal -> CReal
realExp x = makeCRealFun2 (expUniformCts (1+intApprox x)) x

{-Requires that abs(a!!i+1) < abs(a!!i) and the sign of the terms alternate -}
alternatingSeries :: [Base] -> Complete Base
alternatingSeries a eps = sum partSeries
 where
  partSeries = (takeWhile (\x -> (abs x) > eps) a)

rationalSin :: Base -> Base -> CReal
rationalSin tol x | tol <= (abs x) =
                      realBasePolynomial [0, 3, 0, (-4)] (rationalSin tol (x/3))
                  | otherwise = CReal (alternatingSeries series) 0
 where
  series = fst $ unzip $ iterate (\(t,n) -> (-t*(x^2)/(n^2+n),n+2)) (x, 2)

sinCts :: UniformCts Base (Complete Base)
sinCts = UniformCts id (approx . rationalSin radius)
```

```
realSlowSin :: CReal -> CReal
realSlowSin = makeCRealFun2 sinCts

realSin :: CReal -> CReal
realSin x | 0==m = realSlowSin x'
          | 1==m = realSlowCos x'
          | 2==m = negate $ realSlowSin x'
          | 3==m = negate $ realSlowCos x'
 where
  n = intApprox (x / realPi2)
  m = n `mod` 4
  x' = x - (realScale (fromInteger n) realPi2)

rationalCos :: Base -> Base -> CReal
rationalCos tol x = realBasePolynomial [1, 0, (-2)] (rationalSin tol (x/2))

cosCts :: UniformCts Base (Complete Base)
cosCts = UniformCts id (approx . rationalCos radius)

realSlowCos :: CReal -> CReal
realSlowCos = makeCRealFun2 cosCts

realCos :: CReal -> CReal
realCos x | 3==m = realSlowSin x'
          | 0==m = realSlowCos x'
          | 1==m = negate $ realSlowSin x'
          | 2==m = negate $ realSlowCos x'
 where
  n = intApprox (x / realPi2)
  m = n `mod` 4
  x' = x - (realScale (fromInteger n) realPi2)

{- computes ln(x).  only valid for 1<=x<2 -}
rationalSmallLn :: Base -> CReal
rationalSmallLn x = assert (1<=x && x<=(3/2)) $
  makeCReal $
  alternatingSeries (zipWith (*) (poly 1) (tail (powers (x-1))))
 where
  poly n = (1/n):(-1/(n+1)):(poly (n+2))

{- requires that 0<=x -}
rationalLn :: Base -> CReal
rationalLn x | x<1 = negate (posLn (recip x))
             | otherwise = posLn x
 where
  ln43 = rationalSmallLn (4/3)
  ln2 = wideLn 2
  {- good for 1<=x<=2 -}
  wideLn x | x < (3/2) = rationalSmallLn x
           | otherwise = (rationalSmallLn ((3/4)*x)) + ln43
  {- requires that 1<=x -}
  posLn x | n==0 = wideLn x'
          | otherwise = (wideLn x') + (realScale n ln2)
    where
     (x',n) = until (\(x,n) -> (x<=2)) (\(x,n) -> (x/2,n+1)) (x,0)

{- domain is [nonZero, inf) -}
lnUniformCts :: Base -> UniformCts Base (Complete Base)
lnUniformCts nonZero = UniformCts mu f
 where
  f x = approx $ rationalLn (max x nonZero)
  mu eps = assert (nonZero > 0) $ eps*nonZero

realLnWitness :: CReal -> Base -> CReal
realLnWitness x nonZero = makeCRealFun2 (lnUniformCts nonZero) x

realLn :: CReal -> CReal
realLn = realLnWitness x (proveNonZero x)

{- only valid for (abs x) < 1 -}
rationalSmallArcTan :: Base -> CReal
rationalSmallArcTan x = assert ((abs x)<(1/2)) $ makeCReal $
  alternatingSeries (zipWith (\x y->x*(y^2)) (series 0) (powers x))
 where
  series n = (x/(n+1)):(-x/(n+3)):(series (n+4))

rationalArcTan :: Base -> CReal
rationalArcTan x | x <= (-1/2) = negate $ posArcTan $ negate x
                 | otherwise = posArcTan x
 where
  {-requires (-1/2) < x-}
  posArcTan x | 2 < x = realPi2 - rationalSmallArcTan (recip x)
              | (1/2) <= x = realPi4 + rationalSmallArcTan y
              | otherwise = rationalSmallArcTan x
    where
     y = (x-1)/(x+1)

arcTanCts :: UniformCts Base (Complete Base)
arcTanCts = UniformCts id (approx . rationalArcTan)

realArcTan :: CReal -> CReal
realArcTan = makeCRealFun2 arcTanCts
```

```
{- Computes x * Pi -}
{- http://aemes.mae.ufl.edu/~uhk/PI.html -}
scalePi :: Base -> CReal
scalePi x =
  ((realScale (x*48) (rationalSmallArcTan (1/38))) +
   (realScale (x*80) (rationalSmallArcTan (1/57))) +
  ((realScale (x*28) (rationalSmallArcTan (1/239))) +
   (realScale (x*96) (rationalSmallArcTan (1/268)))))

real2Pi = scalePi 2
realPi = scalePi 1
realPi2 = scalePi (1/2)
realPi4 = scalePi (1/4)

nestedBalls :: [(Base, Base)] -> Complete Base
nestedBalls ((lb,ub):bs) eps | (ub - lb) < 2*eps = (ub+lb)/2
                             | otherwise = nestedBalls bs eps

{- My algorithm -}
{-
rationalSqrt :: Base -> CReal
rationalSqrt x = makeCReal (nestedBalls (iterate betterBounds (0,(1+x)/2)))
 where
  betterBounds (lb, ub) = assert (lb <= lb' && lb' <= ub' && ub' <= ub) $
    (lb', ub')
   where
    intercept a b = (x + a*b)/(a+b)
    lb' = intercept lb ub
    ub1 = intercept lb lb
    ub2 = intercept ub ub
    ub' = if (lb > 0) then (min ub1 ub2) else ub2
-}
{- Freek's algorithm -}
rationalSqrt :: Base -> CReal
rationalSqrt n | n < 1 = realScale (1/2) (rationalSqrt (4*n))
               | 4 <= n = realScale 2 (rationalSqrt (n/4))
               | otherwise = makeCReal (\eps -> f eps)
 where
  f eps = fst $ until (\(x,err) -> err <= eps) newton ((1+n)/2, 1/2)
  newton (x,err) = (approxBase x' e1, e1*2)
   where
    x' = (n+x^2)/(2*x)
    e1 = err^2/2

sqrtCts :: UniformCts Base (Complete Base)
sqrtCts = UniformCts (^2) (approx . rationalSqrt)

realSqrt :: CReal -> CReal
realSqrt = makeCRealFun2 sqrtCts

instance Floating CReal where
 exp = realExp
 log = realLn
 pi = realPi
 sin = realSin
 cos = realCos
 atan = realArcTan
 sqrt = realSqrt
 sinh x = realScale (1/2) (exp x - (exp (-x)))
 cosh x = realScale (1/2) (exp x + (exp (-x)))
 asin x = atan (x/sqrt(realTranslate 1 (negate (realPowerInt x 2))))
 acos x = realPi2 - asin x
 acosh x = log (x*sqrt(realTranslate (-1) (realPowerInt x 2)))
 asinh x = log (x*sqrt(realTranslate 1 (realPowerInt x 2)))
 atanh x = realScale (1/2)
   (log ((realTranslate 1 x) / (realTranslate 1 (negate x))))

{- testing stuff is below -}
test0 = makeCReal id

answer n x = shows (intApprox (realScale (10^n) x))
  "x10^-"++(show n)

sumRealList :: [CReal] -> CReal
sumRealList [] = realBase 0
sumRealList l = makeCReal (\eps -> sum (map (\x -> approx x (eps/n)) l))
 where
  n = fromIntegral $ length l
```